\documentclass[10pt]{article}
\usepackage[paper=a4paper,
	top=3cm,
	inner=2.54cm,
	outer=2.54cm,
	bottom=3cm,
	headheight=5ex,
	headsep=5ex]{geometry}
\usepackage[pagewise]{lineno}
\usepackage{amsmath, amssymb}
	\numberwithin{equation}{section}
\usepackage{tikz-cd}
\usepackage{enumitem} 
\usepackage{mathdots}
\usepackage{natbib}
	\setcitestyle{numbers,square}
\usepackage[bookmarks,colorlinks ]{hyperref}
	\hypersetup{linkcolor =blue,
	filecolor=red}

\usepackage{amsthm}
\newtheorem{thm}{Theorem}[section]
\newtheorem{prop}[thm]{Proposition}
\newtheorem{lem}[thm]{Lemma}
\newtheorem{cor}[thm]{Corollary}
\newtheorem{fact}[thm]{Fact}

\theoremstyle{definition}

\theoremstyle{remark}
\newtheorem{rmk}{Remark}[subsection]
\newtheorem{eg}{Example}

\newenvironment{pf}[1][]{{\noindent\it{Proof#1}}. \small}{\hfill \qed}

\newcommand{\Z}{\mathbb Z}
\newcommand{\R}{{\mathbb R}}
\newcommand{\C}{{\mathbb C}}
\newcommand{\g}{{\mathfrak g}}
\renewcommand{\k}{{\mathfrak k}}
\newcommand{\h}{{\mathfrak h}}
\renewcommand{\t}{{\mathfrak t}}
\renewcommand{\a}{{\mathfrak a}}
\renewcommand{\H}{{\mathbb H}}
\newcommand{\Mod}{{\textrm{-}\mathrm{Mod}}}
\newcommand{\m}{{\mathfrak m}}
\newcommand{\For}{\operatorname{For}}
\newcommand{\Jac}{\operatorname{Jac}}
\newcommand{\St}{{\mathrm{St}}}
\newcommand{\la}{\left\langle}
\newcommand{\ra}{\right\rangle}

\newcommand{\pk}{k}

\begin{document}
\title{Comparing Translation and Jaquet functor over general linear groups}
\author{Chang Huang}
\maketitle
\begin{abstract}
	Kei Yuen Chan and Kayue Daniel Wong constructed a functor  from the category of Harish-Chandra modules of $\mathrm{GL}(n, \C)$ to the category of modules over graded Hecke algebra $\mathbb H_m$ of type A.
	This functor has several nice properties, such as compatible with parabolic inductions, and preserving standard and irreducible objects.
	Based on their results, we show this functor relates translation functor on the real side and Jacquet functor on the $p$-adic side.
\end{abstract}


\newcommand{\Q}{\mathbb Q}
\newcommand{\Irr}{\mathrm{Irr}}
\newcommand{\Hom}{\mathrm{Hom}}
\newcommand{\N}{\mathbb N}
\newcommand{\HC}{{\mathcal{HC}}}
\newcommand{\GL}{{\mathrm {GL}}}
\section{Introduction}
	A recent work \cite{Chan23} studies the Lefschetz principle between $\GL(n, \C)$ and $\GL(m, \Q_p)$.
	It constructs a family of functors $\Gamma_{n, m}$ from representations of $\GL(n, \C)$ to modules of graded Hecke algebra $\H_m$, and 
	verifies their compatibility with many classical operators parallelly defined on both sides.
	For example, those $\Gamma_{n, m}$ preserve parabolically induced modules, standard modules, irreducible modules, unitary modules and Dirac series.
	
	Another main theorem in \cite{Chan23} is relating the generalized Bernstein-Zelevinsky functors on $p$-adic side with certain tensoring functors on real side.
	This article provides a refinement in a special case: 
	relating the translation functor on real side with the Jacquet functor on $p$-adic side.
	Let's review some definitions, and state our main result as Theorem \ref{main}.

\subsection{Translation functor}
	Let $G$ be a real Lie group, $\g$ be its complex Lie algebra, $U(\g)$ be the universal enveloping algebra of $\g$, and $Z(\g)$ be the center of $U(\g)$.
	Take a multiplicative functional $\chi$ over $Z(\g)$, then
	there is an Abelian category $\HC_\chi(G)$ of Harish-Chandra modules over $G$ with generalized infinitesimal character $\chi$.
	Fix a Cartan subalgebra $\h$ of $\g$, then $\chi$ determines uniquely a Weyl group orbit in $\h^*$ according to Harish-Chandra isomorphism $Z(\g) \cong \C[\h]^W$.
	Conversely, each linear functional $\lambda \in \h^*$ determines a multiplicative functional $\chi_\lambda$ over $Z(\g)$, 
	and we will abbretive $\HC_{\chi_\lambda}(G)$ by $\HC_\lambda(G)$.
	
	Let $\Lambda \in \h^*$ be an integral weight, and $\pi_\Lambda$ be the finite dimensional representation of $G$ with extreme weight $\Lambda$.
	There defines the translation functor 
	$T_\lambda^{\lambda+ \Lambda}: \HC_\lambda(G) \to \HC_{\lambda+\Lambda}(G)$, $X \mapsto (\lambda + \Lambda)$-primary subspace of $X \otimes_\C \pi_\Lambda$ with respect to $Z(\g)$.
	We refer to \cite[Chapter VII]{KV95} for detailed properties of the translation functor.

\subsection{Jacquet functor}\label{subsec: H_m}
	We define the graded Hecke algebra $\H_m$ of type A as follows.
	It is an algebra generated by polynomial ring $\C[y_1, \cdots, y_m]$ and group algebra $\C[S_n]$, and subject to the relations:
	$s_i y_i - y_{i+1} s_i =1$,
	$s_i y_j - y_j s_i =0, j \not= i, i+1$.
	Here $s_i= (i, i+1)$ are transpositions that form the Coxeter generators of $S_n$.
	Denote the Abelian category of finite dimensional modules over $\H_m$ by $\H_m\Mod$.
	
	There is an embedding of algebra $\H_{m-1} \to \H_m$ given by $y_i \mapsto y_{i+1}$, $s_i \mapsto s_{i+1}$.
	Its image is commutative with $y_1 \in \H_m$.
	Then for $a \in \C$, there defines the functor $\Jac_a: \H_m\Mod \to \H_{m-1}\Mod$ by taking the generalized $a$-eigenspace with respect to $y_1$.
	Our definition is an analogue of the operators $e_i^*$ defined for affine Hecke algebra in the front of \cite[Section 8]{Gro99}.

\subsection{Main results}\label{result-and-K}
	For $G = \GL(n, \C)$, we denote the Abelian category of Harish-Chandra modules over $G$ by $\HC_n$, and 
	abbreviate $\HC_\lambda(G)$ to $\HC_\lambda$.
	
	Let $H \subseteq G$ be the Cartan subgroup consisting of diagonal matrices.
	It has a natural decomposition $H \cong T \times A$ into compact part $T$ and split part $A$.
	On the level of complex Lie algebras we also have $\h \cong \t \oplus \a$,
	so a weight $\lambda \in \h^*$ can be writen as $(\mu, \nu) \in \t^* \oplus \a^*$.
	Since $\t \cong \C^n \cong \a$ naturally, $\mu$ and $\nu$ can be writen further in to $n$-tuples of complex numbers.
	Let $e_\pk \in \C^n$ denote the $\pk$-th unit vecor whose $\pk$-th entry is $1$ and others are $0$.
\begin{thm}\label{main}
	The following diagrams commute:
	\begin{equation}\label{diagram: compability}
	\begin{tikzcd}
	{\HC_{(\mu, \nu)}} \arrow[r, "{\Gamma_{n, m}}"] \arrow[d, "{T_{(\mu, \nu)}^{(\mu - e_\pk, \nu + e_\pk)}}"'] & 
	\H_m \Mod \arrow[d, "{\Jac_{\frac {\nu_\pk - \mu_\pk +1} 2}}"] \\
	{\HC_{(\mu - e_\pk, \nu + e_\pk)}} \arrow[r, "{\Gamma_{n, m-1}}"]                                                                    &
	\H_{m-1} \Mod                                             
	\end{tikzcd}.
	\end{equation}
	More precisely, there is a natural isomorphism
	$\Gamma_{n, m-1} \circ T_{(\mu, \nu)}^{(\mu - e_\pk, \nu + e_\pk)}
	\simeq
	\Jac_{\frac {\nu_\pk - \mu_\pk +1} 2} \circ \Gamma_{n, m}$.
	
	Dually, there is also a natural isomorphism
	$\Gamma_{n, m-1} \circ T_{(\mu, \nu)}^{(\mu - e_\pk, \nu - e_\pk)} \simeq
	\Jac^{\frac {\nu_\pk + \mu_\pk +1} 2} \circ \Gamma_{n, m}.$
	Here, $\Jac^a : \H_m\Mod \to \H_{m-1}\Mod$ denotes the functor that takes the generalized $a$-eigenspace with respect to $y_m\in \H_m$, 
	and regards it as an $\H_{m-1}\Mod$ via the natural embedding $\H_{m-1} \to \H_m$, $y_i \mapsto y_i$, $s_i \mapsto s_i$.
\end{thm}

\subsection{Relation to other results}
	\cite[Theorem 1.2]{Chan23} asserts the natrual isomorphism
	$\mathbf{BZ}_\tau \circ \Gamma_{n, m} (-) \simeq \Gamma_{n, m-k} ( - \otimes \mathbb S_\tau(V))$, where $\tau \in  \Irr(S_k)$, and
	\begin{itemize}
	\item $\mathbf{BZ}_\tau (-) = \Hom_{S_k}( \tau, -): \H_m \Mod \to \H_{m-k} \Mod$, 
	\item $\mathbb S_\tau(V)$ is a finite-dimensional irreducible representation of $G$.
	\end{itemize}
	Set $\tau=$ trivial representations of $S_1$, then $\mathbf{BZ}_\tau$ is the forgetful functor, and $\mathbb S_\tau(V) = V$.
	It's easy to check that the translation functor $T_{(\mu, \nu)}^{(\mu - e_\pk, \nu+ e_\pk)}$ is a direct summand of  the tensoring functor $- \otimes V$, and
	the Jaqcuet functor $\Jac_{\frac{\nu_\pk - \mu_\pk+1} 2 }$ in is a direct summand of the forgetful functor.
	Thus, our main Theorem is a refinement of theirs.
	We will also mention this after Fact \ref{BZ-and-tensor}.
	
	There is a geometric version on the Langlands dual side given in \cite{25}.
	The authors show that
	\begin{itemize}
	\item (Theorem 1.1.2 and 3.3.13) $\Gamma_{n, m}$ is dual to the pull back $\zeta^*$ by a locally closed embedding,
	\item (Theorem 2.2.1) translation functor is dual to some push-pull functor, and
	\item (Theorem 2.4.8) Jacquet functor is dual to the Lusztig's inductioin.
	\end{itemize}
	Moreover, their Theorem 1.1.1 and 4.3.1 asserts that $\zeta^*$ intertwines the push-pull functor and Lusztig's induction.
	This result is indeed dual to our Theorem \ref{main}.

\subsection{Outline}
	The main body of this article will be for the commutativity of \eqref{diagram: compability}.
	Section \ref{sec: representations} collects several necessary definitions and results in \cite{Chan23}.
	Section \ref{sec: Nat} provides the natural transformation for \eqref{diagram: compability} in Proposition \ref{nat-iso},
	and reduces the verification to case of principle series representation.
	Section \ref{sec: eigenvalue} contains the most technical calculations.
	Section \ref{dual} discusses the dual diagram of \eqref{diagram: compability}.

\subsection*{Acknowledgement}
	The author would like to thank Professor Bin Xu, Taiwang Deng and Qixian Zhao for patient guidance.

\renewcommand{\d}{\, \mathrm d}
\renewcommand{\Re}{\operatorname{Re}}
\newcommand{\gl}{{\mathfrak {gl}}}
\newcommand{\U}{{\mathrm U}}
\newcommand{\Ind}{\operatorname{Ind}}
\newcommand{\sgn}{\operatorname{sgn}}
\newcommand{\Ad}{\operatorname{Ad}}
\newcommand{\Tr}{\operatorname{Tr}}
\newcommand{\End}{\operatorname{End}}
\newcommand{\Id}{\operatorname{Id}}
\section{Preliminaries on Representations}\label{sec: representations}
	In this section we review the basic representation theory for $\GL(n, \C)$ and $\H_m$.
	The subsequent materials are mainly copied from \cite{Chan23}.

\subsection{Representations of $\GL(n, \C)$}
	Although $G=\GL(n, \C)$ has complex structure, we understand it as a real Lie group, 
	and study its representations through $(\g, K)$-modules.

\subsubsection{Complex Lie algebra}
	Denote the Lie algebra of $G$ by $\g_0 = \gl(n, \C)$,
	and its complexfication by $\g = \g_0 \otimes_\R \C$.
	There is a natural decomposition of real vector space $\g = \g_0 \oplus j \g_0$, 
	where $j$ denotes the action of $\sqrt{-1}\in \C$ on $\g$.
	However, it is not a decomposition of complex Lie algebra.
	Define $\phi^L, \phi^R: \g_0 \to \g$ by
	\[\phi^L(E) = \frac 1 2 (E- ji E),\quad
	\phi^R(E)= \frac1 2 (\bar E+ ji \bar E),\]
	where $i$ denotes the multiplication of $\sqrt{-1}\in \C$ on $\g_0$.
	Then one can check by direct computation that
	\begin{itemize}
	\item $\phi^L$, $\phi^R$ are $\C$-linear maps,
	\item $\phi^L$, $\phi^R$ preserve Lie brackets, and
	\item $(\phi^L, \phi^R): \g_0 \times \g_0 \to \g$, $(E, E') \to \phi^L(E) + \phi^R(E')$ gives an isomorphism of complex Lie algebra.
	\end{itemize}
	Denote the element $\phi^L(E) + \phi^R(E') \in \g$ by $(E, E') \in \g_0 \times \g_0$.
	Under this coordinate, 
	\begin{itemize}
	\item the natural embedding $\g_0 \to \g_0 \otimes_\R \C \cong \g_0 \times \g_0$ becomes $E \mapsto (E, \bar E)$.
	\item the conjugation of $\g = \g_0 \otimes_\R \C$ is $(E, E') \mapsto (\bar E', \bar E)$, and
	\item the Cartan involution $\theta$ of $\g$ is $(E, E') \mapsto (-E'^t, -E^t)$.
	\end{itemize}
	The maximal compact subgroup of $G$ is $K = \U(n)$;
	its real Lie algebra $\k_0$ coincide with the $\theta$-fixed points in $\g_0$.
	
	Recall that $H \subseteq G$ consists of diagonal matrices, and has compact part $T$, split part $A$.
	The complex Lie algebra $\h \subseteq \g$ of $H$ is sent to $\h_0 \times \h_0$ under $\g \cong \g_0 \times  \g_0$.
	In the mean time, the complex Lie algebras $\t, \a$ of $T$, $A$ are isomorphic to $\h_0 \cong \C^n$ via
	\[\h_0  \to \t, H \mapsto (H, -H), \quad
	\h_0 \to \a, H \mapsto (H, H).\]
	A weight $\lambda \in \h^*$ can be denoted by $(\lambda_L, \lambda_R)$ under coordinate $\h \cong \h_0 \oplus \h_0$, and by $(\mu, \nu)$ under $\h \cong \t \oplus \a$.
	Identify $\t$, $\a$ with $\h_0$ as above,
	then the coordinates $(\lambda_L, \lambda_R)$ and $(\mu, \nu)$ are related by
	$\mu= \lambda_L - \lambda_R$,
	$\nu= \lambda_L + \lambda_R$.
	We will switch between these two coordinates freely.
	In particular, the translation functor $T_{(\mu, \nu)}^{(\mu-e_\pk, \nu+ e_\pk)}$, Jacquet functor $\Jac_{\frac{\nu_\pk - \mu_\pk +1} 2}$ in diagram \eqref{diagram: compability} can also be denoted by $T_{(\lambda_L, \lambda_R)}^{(\lambda_L, \lambda_R+e_\pk)}$, $\Jac_{\lambda_{R, \pk} + \frac 1 2}$.

\subsubsection{Principle series representation}\label{subsubsec: Principle series}
	Let $B \subseteq G$ consist of upper triangular matrices; it is a Borel subgroup. 
	If $\mu \in \t^*$ is integral, then $\lambda = (\mu, \nu) \in \t^* \oplus \a^* \cong  \h^*$ can be lifted to a character of $H = TA$, and 
	further to a one dimensional representation $\C_\lambda$ of $B$.
\begin{eg}
	The modular character
	\[\delta: B \to \C^\times,
	\begin{pmatrix}
		b_1	&	&*	\\
			&\ddots&	\\
			&	&b_n
	\end{pmatrix}
	\mapsto
	\prod_{i=1}^n |b_i|^{2i-1-n}\]
	can be obtained from $(0, \cdots, 0, n-1, \cdots, 1-n) \in \C^n \oplus \C^n \cong \t^* \oplus \a^* \cong \h^*$ in this way.
	Under the coordinate $\C^n \oplus \C^n \cong \h_0^* \oplus \h_0^* \cong \h^*$, it is also represented by $(\delta_L, \delta_R)$ with
	$\delta_L = \delta_R = (\frac{n-1}2, \cdots, -\frac{n-1}2)$.
\end{eg}
	For such a one dimensional representation $\C_\lambda$, 
	let $\Ind_B^G \C_\lambda $ denote the completion of following function space 
	\[\{ f\in C^\infty(G) \mid f(bg) = \delta(b)^{\frac 1 2} \lambda(b) f(g),
	\forall g \in G, b \in B\}\] 
	under inner product
	\[\la f_1, f_2 \ra = \int_K \overline{f_1(k)} f_2(k) \d k ,\]
	and $G$ act by right translation.
	It is called a principle series representation, and denoted by $X(\lambda)$ in \cite{Chan23}.
	By taking $K$-finite subspace, we can also obtain its $(\g, K)$-module version.
	
	These principle series $X(\lambda) = X(\lambda_L, \lambda_R)$ play an important role in the classfication of irreducible representations.
	Each $X(\lambda)$ has a unique subquotient $\bar X(\lambda) = \bar X(\lambda_L, \lambda_R)$ that contains the $K$-type with extremal weight $\mu = \lambda_L - \lambda_R$ (\cite{PRV67});
	we call $\bar X(\lambda)$ the Langlands subquotient of $X(\lambda)$.
\begin{fact}[\cite{Zhe74}]
	For a weight $\lambda \in \h^*$, we will use both the coordinates $(\lambda_L, \lambda_R) \in \h_0^* \oplus \h_0^*$ and $(\mu, \nu) \in \t^* \oplus \a^*$.
	\begin{enumerate}[label=$(\arabic*)$]
	\item The (isomorphsm class of) irreducible Harish-Chandra modules of $G$ is exhausted by $\bar X(\mu, \nu)$, 
	where $\mu \in \t^*$ is integral.
	\item Two modules $\bar X(\lambda_L, \lambda_R)$ and $\bar X(\lambda_L', \lambda_R')$ are isomorphic if and only if 
	there is $w \in W(\g_0, \h_0) \cong S_n$ such that $w \lambda_L = \lambda_L',$ $w\lambda_R = \lambda_R'$.
	\item For $w\in W(\g_0, \h_0)$, $X(w\lambda_L, w\lambda_R)$ has the same composition series with $X(\lambda_L, \lambda_R)$.
	\item The representation $X(\mu, \nu)$ is tempered if and only if $\Re \nu =0$, and in this case $X(\mu, \nu) = \bar X(\mu, \nu)$.
	\end{enumerate}
\end{fact}
	The following facts are also standard, see for example Proposition 11.43 and Theorem 11.198 in \cite{KV95}:
	\begin{itemize}
	\item $X(\lambda)$ has infinitesimal character represented by $\lambda \in \h^*$;
	\item if $X(\lambda)$ is standard, which means $\Re \nu \in \a_0^* \cong \R^n$ is decreasing,
	then $\bar X(\lambda)$ is the unique quotient of $X(\lambda)$.
	\end{itemize}
	
	Fix the weight $\lambda \in \h^*$.
	Recall that $\HC_\lambda$ is the Abelian category of Harish-Chandra modules over $\GL(n, \C)$ with generalized infinitesimal character (represented by) $\lambda$.
	The above facts imply that irreducible objects in $\HC_\lambda$ can be exhausted by quotients of standard principle series.

\subsection{Representations of $\H_m$}\label{rep-of-H}
	There is a natural embedding $\H_{m_1} \otimes_\C \H_{m_2} \to \H_{m_1+m_2}$.
	Define the parabolic induction 
	$\H_{m_1}\Mod \times \H_{m_2}\Mod \to \H_{m_1+m_2}\Mod$ 
	via 
	\[M_1 \times M_2: =
	 \H_{m_1+m_2} \otimes_{\H_{m_1} \otimes_\C \H_{m_2}} 
	 ( M_1 \boxtimes M_2).\]
	
	Irreducible representations of $\H_1= \C[y]$ are all one-dimentional character, 
	and exhausted by the evaluation homomorphism $\psi_c$, $c\in \C$ with $\psi_c(y) = c$.
	Irreducible representations of $\H_m$ ($m \geqslant 2$)  are classified by multi-segements (multi-sets of segements) in $\C$.
	By a segement $[a, b]$ with $a, b\in \C$, $b-a\in \Z$ we mean the set $\{a, a+1, \cdots b\}$;
	if $b<a$, then this set is empty.
	For a segement $\Delta = [a, b]$ with length $m = b-a +1 \geqslant 1$, there is a representation of $\H_m$ given by
	$\psi_a \times \cdots \times \psi_b$.
	It has a unique irreducible quotient denoted by  $\St(\Delta)$.
	For $\Delta= \varnothing$, we define $\St(\Delta)$ to be the trivial representation $\C$ of $\H_0 = \C$.
	
\begin{lem}\label{Leb-of-Jac}
	Let $\Delta = [a, b] \subset \C$ be a segement,
	then representation $\St(\Delta)$ is one-dimensional, with $S_m$ acting by $\sgn: w \mapsto (-1)^{l(w)}$, and $y_1, \cdots, y_m$ acting by $a, \cdots, b$.
	In particular,
	\[\Jac_x (\St(\Delta))= \left\{\begin{aligned}
	&\St({}^-\Delta),	&\textrm{if } x=a,\\
	&0,	&\textrm{otherwise},
	\end{aligned}\right.\]
	where ${}^-\Delta = [a+1, b]$ is obtained from $\Delta$ by deleting its beginning.
\end{lem}
\begin{pf}\begin{subequations}
	The case $\Delta= \varnothing$ and $\{a\}$ is trivial.
	For other $\Delta$, prove by induction on the length $m =b -a +1 \geqslant 2$.
	Let $\Delta^- = [a, b-1]$, then $\St(\Delta)$ is the irreducible quotient of $\St(\Delta^-) \times \psi_b$.
	Under the induction bypothesis, we know $\St(\Delta^-) \times \psi_b = \Ind_{S_{m-1}}^{S_m}(\sgn)$ as a representation of $S_m$.
	Let's take the set $\{s_{j} \cdots s_{m-1} \mid j = 1, \cdots, m\}$ as representatives for $S_m / S_{m-1}$, then 
	the representation space for $\Ind_{S_{m-1}}^{S_m}(\sgn)$ can be taken to be
	\[\bigoplus_{j=1}^m \C [ s_{j} \cdots s_{m-1} ].\]
	It has a one-dimensional $S_m$-subrepresentation, spanned by
	$e = \sum_j  (-1)^{j-1} [ s_j \cdots s_{m-1}]$;
	its orthogonal complement (under the standard inner product) $E$ is also invariant under $S_m$.
	
	One calculates from definition that
	\begin{align}\label{eq: action-of-y-on-St-1}
	y_1 \cdot [s_1 \cdots s_{m-1}] &= 
	\sum_{i=1}^{m-1}(-1)^{i-1} [ s_{i+1} \cdots s_{m-1}]
		+ b [s_1 \cdots s_{m-1} ], \\
	\label{eq: action-of-y-on-St-2}
	y_1 \cdot [s_j \cdots s_{m-1} ] &=a [s_j \cdots s_{m-1}], \quad j\geqslant 2.
	\end{align}
	Since $E$ consists of the vectors $\sum_j c_j [s_j \cdots s_{m-1}]$ with $\sum_j (-1)^{j-1} c_j = 0,$
	one check from the above formula that this space is invariant under $y_1$.
	Moreover, from the generating relation $y_{i+1} = s_i y_i s_i - s_i$, one can deduce by induction that $E$ is invariant under all $y_i \in \H_m$;
	thus, $E$ is a $\H_m$-submodule of $\St(\Delta^-) \times \psi_b$.

	According to the definition, $\St(\Delta)$ is the quotient of $\St(\Delta^-) \times \psi_b$ by $E$.
	In particular, it is one dimensional, and can be spanned by the image of $e$.
	The group $S_m \subset \H_m$ acts on it by $\sgn$, since each $s_i \in S_m$ acts by $-1$.
	From the formulas \eqref{eq: action-of-y-on-St-1} and \eqref{eq: action-of-y-on-St-2}, one calculates that
	$y_1 \cdot e - a e\in E,$
	so $y_1$ acts on $\St(\Delta)$ by the scalar $a$.
	Then, from the generating relation $y_{i+1} = s_i y_i s_i - s_i$, one deduces by induction that $y_{i+1}$ acts on $\St(\Delta)$ by the scalar $a+ i$.
%
\end{subequations}\end{pf}

	Two segements $\Delta_1, \Delta$ are said to be linked, if $\Delta_1 \cup \Delta_2$ is still a segement, and they are not in containment relation.
	If $\Delta_1 = [ a_1, b_1]$ and $\Delta_2 = [a_2, b_2]$ are linked, and $a_1 < a_2$, then write $\Delta_1 < \Delta_2$.
	For a multi-segement $\m$ in $\C$, enumerate its segements by $\Delta_1, \cdots, \Delta_n$ properly, such that they meet the requirement:
	if $\Delta_i < \Delta_j$, then $i>j$.
	Then $\St(\Delta_1) \times \cdots \times \St(\Delta_n)$ has a unique irreducible qoutient, and
	we denote it by $\mathrm{Speh}(\m)$.
	Such representations exhaust all irreducible representations of $\H_m$;
	this is a standard fact.
	Take the graded ring
	\[K\H = \bigoplus_{m \geqslant 0} K\H_m,\]
	with multiplication given by parabolic induction.
	It is a polynomial ring (over $\Z$), with indeterminates given by $\{\St(\Delta) \mid \Delta \subset \C$ segements$\}$.
	The Jacquet functor $\Jac_a$ induces a derivation over the graded ring $K\H$;
	more explicitly, $\forall [M_1], [M_2] \in K\H$, $\Jac_a( [M_1] \times [ M_2]) = (\Jac_a[ M_1]) \times [M_2] + [M_1] \times (\Jac_a[M_2])$.
	However, the above facts in this paragraph will not be used in this article.

\subsection{The functor $\Gamma_{n, m}$}\label{Def-of-Gamma}
	Let $V$ denote the conjugate standard representation of $G= \GL(n, \C)$.
	This means $V= \C^n$, and $g\in G$ acts on it by left multiplication of the matrix $\bar g$.
	
	For $X\in \HC_\lambda$, take $\Gamma_{n, m}(X) = (X \otimes V^{\otimes m})^K$.
	To make it into an $\H_m$-module, \cite{Chan23} defines an action of $\H_m$ on $X\otimes V^{\otimes m}$, which is commutative with the action of $G$,
	thus preserves its $K$-fixed points.
	Here we present the conceret actions of $s_k, y_l \in \H_m$ on $v_0 \otimes v_1 \otimes \cdots \otimes v_m \in X \otimes V^{\otimes m}$ as follows: 
	let $E_{i, j} \in \g_0$ denote the $n \times n$ matrix 
	whose only non-zero entry is $1$ at $(i, j)$, then
	\begin{subequations}\label{action-of-H}
	\begin{align}
	s_k \cdot v_0 \otimes \cdots \otimes v_m =&
	-\sum_{1\leqslant i , j \leqslant n} 
		v_0 \otimes \cdots \otimes (0, E_{i,j}) v_k 
		\otimes (0, E_{j, i}) v_{k+1} \otimes \cdots \otimes v_m,
	\label{eq: action-of-s}	\\
	y_l \cdot v_0 \otimes \cdots \otimes v_m =&
	\sum_{0 \leqslant x <l} \sum_{1\leqslant i , j \leqslant n}
		v_0 \otimes \cdots \otimes (0, E_{i,j}) v_x \otimes \cdots 
		\otimes (0, E_{j, i}) v_l \otimes \cdots \otimes v_m
		\notag	\\
		&+ \frac n 2 v_0 \otimes \cdots \otimes v_m.	
	\label{eq: action-of-y}
	\end{align}
	\end{subequations}
	
\begin{fact}\label{fact: action-of-H}
	The above action of $s_k$, $y_l$ on $X \otimes V^{\otimes m}$
	satisfies the generating relation of $\H_m$.
	This action commutes with $G$.
	Moreover, the action of $s_k$on $X\otimes  V^{\otimes m}$ can be expressed as
	\begin{equation}\label{rmk: action-of-s}
	v_0\otimes \cdots \otimes v_k \otimes v_{k+1} \otimes \cdots v_m
	\mapsto
	- v_0\otimes \cdots \otimes v_{k+1} \otimes v_k \otimes \cdots v_m.
	\end{equation}
\end{fact}
\begin{pf}
	For the well-defineness, Chan-Wong refers to \cite[Lemma 2.2.1]{AS98}.
	One can also check it from \cite[Lemma 2.3.2]{CT11} together with the explict formula for the action of $s_k$.
	
	The commutativity is explained in the end of \cite[Subsection 3.2]{Chan23}.
	The main point is that $\forall g\in G$, the acton of $\Ad(g)$ on $\g_0 \otimes \g_0$ stablizes
	\[\Omega := \sum_{1\leqslant i, j \leqslant n} E_{i, j} \otimes E_{j, i}.\]
	It can be deduced as follows.
	The trace paring $\g_0 \times_\C \g_0 \to \C$, $(E, E') \to \Tr(E E')$ is invariant under $G$, so 
	it induces a $G$-equivariant isomorphism $\g_0 \otimes_\C \g_0 \to \g_0^* \otimes_\C \g_0 = \End_\C(\g_0)$,
	where $g \in G$ acts on $\g_0 \otimes_\C \g_0$ by $E \otimes E' \mapsto \Ad(g) E \otimes \Ad(g) E'$, and 
	acts on $\End_\C(\g_0)$ by $\phi \mapsto \Ad(g) \circ \phi \circ \Ad(g)^{-1}$.
	The basis $\{E_{i, j}\}$ is dual to $\{E_{j, i}\}$ under this pairing, so 
	$\Omega$ is sent to $\Id_{\g_0}$ under $\g_0 \otimes_\C \g_0 \to \End_\C(\g_0)$. 
	Now $\Id_{\g_0}$ is clearly stablized by every $g\in G$, so is $\Omega$.

	The explict action of $s_k$ is given in \cite[Lemma 2.3.3]{CT11}.
	Here we can check it directly, and 
	it sufficies to show $\phi^R(\Omega)\in \g \otimes_\C \g$ acts on $V\otimes V$ by $x\otimes y \mapsto y \otimes x$. 
	Note that for $E \in \g_0$, $\phi^R(E) = \frac 1 2 ( \bar E + j i \bar E) \in \g = \g_0 \otimes_\R\C$ acts on $V = \C^n$ by left multiplication of the matrix
	$\frac 1 2 ( E + i \cdot \overline{i E} )= E$;
	this is also mentioned in \cite[Lemma 3.1]{Chan23}.
	Take the standard basis $\{e_1, \cdots, e_n\}$ of $V = \C^n$, then
	\[\phi^R(\Omega) (e_k \otimes e_l) = \sum_{1\leqslant i, j \leqslant n} E_{i,j } e_k \otimes E_{j, i} e_l.\]
	The only non-zero summation is $e_l \otimes e_k$, with index $(i, j) = (l, k)$.
	We thus conclude $\phi^R(\Omega)( x\otimes y) = y \otimes x$.
\end{pf}
	
\begin{fact}\label{Std-under-Gamma}
	Let $\lambda= (\mu, \nu) \in \t^* \oplus \a^* \cong \h^*$ be a weight with $\mu \in \t^* \cong \C^n$ integral.
	If all $\mu_i \geqslant 0$, and $\mu_1 +\cdots +\mu_n = m$, 
	then $X(\lambda)$ is sent to 
	$\St(\Delta_1) \times \cdots \times \St(\Delta_n)$
	under the functor $\Gamma_{n, m}$,
	where $\Delta_i = [ \frac {\nu_i - \mu_i + 1} 2, \frac {\nu_i + \mu_i +1} 2 ]$.
	Otherwise, $X(\lambda)$ is sent to $0$.
\end{fact}
\begin{pf}
	This follows from \cite[Lemma 6.1 and Theorem 6.4]{Chan23}.
\end{pf}

\newcommand{\pj}{i}
\section{Natural transformation}\label{sec: Nat}
	In this section we define the natural transformation 
	$\Gamma_{n, m-1} \circ T_{(\lambda_L, \lambda_R)}^{(\lambda_L, \lambda_R+ e_\pk)} 
	\Rightarrow
	\Jac_{\lambda_{R, \pk} + \frac 1 2} \circ \Gamma_{n, m}$,
	and discuss how to prove it to be an isomorphism.

\subsection{Decompose into generalized eigenspaces}	
	Denote the forgetful functor by $\For: \H_m\Mod \to \H_{m-1}\Mod$.
	There is a natural isomorphism 
	\[\For \simeq \bigoplus_{a\in \C} \Jac_a.\]
	On the level of objects, it is merely the decomposition into generalized eigenspaces of $y_1$.

	Consequently, there is a natural isomorphism of functors $\HC_\lambda \to \H_{m-1}\Mod$:
	\[\For \circ \Gamma_{n, m} \simeq \bigoplus_{a\in \C} \Jac_a \circ \Gamma_{n, m}.\]
	One can deduce from Lemma \ref{eigenvalue} that if $a\not \in \{ \lambda_{R, \pj} + \frac 1 2 \mid \pj= 1 , \cdots n \}$, then $\Jac_a \circ \Gamma_{n, m}: \HC_\lambda \to \H_{m-1} \Mod$ is zero.
	Hence, the right hand side of above isomorphism is a finite sum.

\subsection{Consistency of $\H_{m-1}$ actions}
	Recall that $\HC_n$ denotes the Abelian category of Harish-Chandra modules over $G = \GL(n, \C)$.
	There is a natural transformation of functors $\HC_\lambda \to \HC_n$ given by the inclusion
	\[T_{(\lambda_L, \lambda_R)}^{(\lambda_L, \lambda_R+e_\pk)}(X) \hookrightarrow X \otimes V.\]
	Consequently, we have 
	$\Gamma_{n, m-1} \circ T_{(\lambda_L, \lambda_R)}^{(\lambda_L, \lambda_R+ e_\pk)}
	\Rightarrow
	\Gamma_{n, m-1} \circ ( - \otimes V)$.

\begin{fact}\label{BZ-and-tensor}
	As functors $\HC_\lambda \to \H_{m-1}\Mod$,  $\Gamma_{n, m-1} \circ ( - \otimes V)$ coincides with $\For \circ \Gamma_{n, m}$.
	More explictly, $\forall X \in \HC_\lambda$, the action of $\H_{m-1}$ on $\Gamma_{n, m-1} ( X \otimes V) = (X \otimes V^{\otimes m})^K$ coincides with that on $\For\circ \Gamma_{n, m}(X)  = (X \otimes V^{\otimes m})^K$.
\end{fact}
	This fact is a straightforward consequence of \cite[Theorem 8.2]{Chan23}, obtained by setting $\tau$ to be the trivial representation of $S_1$.
	Consequently, our Proposition \ref{nat-iso} can be regarded as a refinement of this established result.
%
%

	Now we obtain the natural transformation of functors $\HC_\lambda \to \H_{m-1}\Mod$
	\[\theta^\pk_a:
	\Gamma_{n, m-1} \circ T_{(\lambda_L, \lambda_R)}^{(\lambda_L, \lambda_R+ e_\pk)}
	\Rightarrow
	\For \circ \Gamma_{n, m}
	\Rightarrow
	\Jac_a \circ \Gamma_{n, m}.\]
	We have to prove 
\begin{prop}\label{nat-iso}
	$\forall X\in \HC_\lambda$, $\theta^\pk_{a, X}:
	\Gamma_{n, m-1} \circ T_{(\lambda_L, \lambda_R)}^{(\lambda_L, \lambda_R+ e_\pk)} (X)
	\to
	\Jac_a \circ \Gamma_{n, m}(X)$ is an isomorphism when $a= \lambda_{R, \pk} + \frac 1 2$, and zero otherwise.
\end{prop}

\subsection{On principle series}
	Here we claim that Proposition \ref{nat-iso} holds for each principle series $X(\lambda)$.
	
\begin{lem}\label{Std-under-Tran}\label{eigenvalue}\label{on-principle}
	For the principle series $X = X(\lambda) = X(\lambda_L, \lambda_R)$, 
	there is a filtration $0 = X_0 \subset X_1 \subset \cdots \subset X_n = X(\lambda)\otimes V$ such that $\forall k = 1, \cdots, n$,
	\[X_k/ X_{k-1} \cong X(\lambda_L, \lambda_R + e_k)\]
	as a representation of $G$.
	Moreover, $y_1 \in \H_m$ preserves each $\Gamma_{n, m-1}(X_k) = (X_k \otimes V^{\otimes (m-1)})^K \subseteq \Gamma_{n, m}(X)$, and 
	acts on $\Gamma_{n, m-1}(X_k/ X_{k-1})$ by generalized eigenvalue $\lambda_{R, k} + \frac 1 2$.
\end{lem}

	The filtration will be given in Subsection \ref{Mackey}, the eigenvalue will be calculated in Section \ref{calculation-over-std}.

\begin{cor}
	For the principle series $X = X(\lambda_L, \lambda_R)$, 
	\[\Gamma_{n, m-1} \circ T_{(\lambda_L, \lambda_R)}^{(\lambda_L, \lambda_R+e_\pk)} (X)=
	\Jac_{\lambda_{R, \pk} + \frac 1 2} \circ \Gamma_{n, m}(X) \]
	as subspaces, and thus $\H_{m-1}$-submodules of $\Gamma_{n, m}(X),$
\end{cor}
\begin{pf}
	Abbreviate $(\lambda_L, \lambda_R + e_\pk)$ to $\lambda'$, and denote $T_\lambda^{\lambda'}\left(X(\lambda)\right)$ by $X'$.
	Recall $X'$ is the $\lambda'$-primary subspace of $X(\lambda) \otimes V$ under $Z(\g)$.
	We have to prove $\Gamma_{n, m-1} (X') = (X' \otimes V^{\otimes (m-1)})^K $ is invariant under action of $y_1$.
	According to Fact \ref{fact: action-of-H}, the endomorphism on $X(\lambda)\otimes V$ given by
	\[x\otimes v \mapsto \sum_{1\leqslant i , j \leqslant n} (0, E_{i, j}) x \otimes (0, E_{j, i}) v\]
	commutes with $G$.
	Then  this endomorphism would preserve the $\lambda'$-primary subspace $X'$ under $Z(\g)$.
	Thus, $y_1\in \H_m$ preserves $X' \otimes V^{\otimes (m-1)} \subset \left(  X(\lambda) \otimes V \right) \otimes V^{\otimes(m-1)}$ due to the definition in \eqref{action-of-H}.
	As a consequence, $y_1$ preserves $\Gamma_{n, m-1}(X') = (X' \otimes V^{\otimes (m-1)})^K 
	\subset (X(\lambda) \otimes V^{\otimes m})^K = \Gamma_{n, m}(X(\lambda))$.
	
	There is an induced filtration $0 = X_0' \subseteq \cdots \subseteq X_n'$ with $X_\pj' = X' \cap X_\pj$.
	Then $X'_\pj / X'_{\pj-1} \hookrightarrow X_\pj/ X_{\pj-1}$ is an embedding onto $\lambda'$-primary subspace.
	On the other hand, $X_\pj/ X_{\pj-1} \cong X(\lambda_L, \lambda_R+ e_\pj)$ is $(\lambda_L, \lambda_R+ e_\pj)$-primary, 
	so $X_\pj'/ X_{\pj-1}'$ is either $0$, or isomorphic to $X(\lambda_L, \lambda_R+ e_\pj)$.
	The latter case occurs if and only if $\lambda' =(\lambda_L, \lambda_R+ e_\pk) $ belongs to the same $S_n\times S_n$-orbit as $(\lambda_L, \lambda_R+ e_\pj)$.
	In particular, if $\lambda_{R, \pj} \not= \lambda_{R, \pk}$, then $X_\pj' / X_{\pj-1}' $ is zero.
	
	Due to Lemma \ref{eigenvalue}, $\Gamma_{n, m-1}(X_\pj)$ is invariant under $y_1$, so is $\Gamma_{n, m-1}(X_\pj') = \Gamma_{n, m-1}(X') \cap \Gamma_{n, m-1}(X_\pj)$;
	moreover, $y_1$ acts on $\Gamma_{n, m-1}(X_\pj'/ X_{\pj-1}') \subseteq \Gamma_{n, m-1} (X_\pj/ X_{\pj-1})$ by generalized eignevalue $\lambda_{R, \pj} + \frac 1 2$.
	However, we have shown that $X_\pj'/ X_{\pj-1}'$ is non-zero only if $\lambda_{R, \pj} = \lambda_{R, \pk}$, 
	so $y_1$ acts on each $\Gamma_{n, m-1}(X_\pj')$ by the unique generalized eigenvalue $\lambda_{R, \pk} + \frac 1 2$.
	Then at least 
	\[\Gamma_{n, m-1}(X') = 
	\Gamma_{n, m-1} \circ T_{(\lambda_L, \lambda_R)}^{(\lambda_L, \lambda_R+e_\pk)} (X) \subseteq
	\Jac_{\lambda_{R, \pk} + \frac 1 2} \circ \Gamma_{n, m}(X).\]
	 The two sides must coincide, since they coincide after taking direct sum with $\lambda_{R, \pk}$ running over $\{ \lambda_{R, \pj} \mid \pj =1, \cdots, n\}$.
\end{pf}

\subsection{On irreducible objects}
	We have explained in Subsubsection \ref{subsubsec: Principle series} that any irreducible object in $\HC_\lambda$ can be realized as a quotient of standard principle series.
	Then Proposition \ref{nat-iso} for irreducible objects can be deduced from the previous conclusion for principle series:
\begin{lem}
	Suppose $X\twoheadrightarrow Y$ is an epimorphism in $\HC_\lambda$, and Proposition \ref{nat-iso} holds for $X$, then it also holds for $Y$.
\end{lem}
\begin{pf}
	For the epimorphism $X \twoheadrightarrow Y$, there is a commutative diagram due to the naturality of $\theta^\pk_a$:
	\[\begin{tikzcd}
	{\Gamma_{n, m-1} \circ T_{(\lambda_L, \lambda_R)}^{(\lambda_L, \lambda_R+e_\pk)}(X)} 
	\arrow[r, "{\theta^\pk_{a, X}}"] 
	\arrow[d, two heads] & 
	{\mathrm{Jac}_a\circ \Gamma_{n, m}(X)} 
	\arrow[d, two heads] \\
	{\Gamma_{n, m-1} \circ T_{(\lambda_L, \lambda_R)}^{(\lambda_L, \lambda_R+e_\pk)}(Y)} 
	\arrow[r, "{\theta^\pk_{a, Y}}"]                  &
	{\mathrm{Jac}_a\circ \Gamma_{n, m}(Y)}                                    
	\end{tikzcd}.\]
	The vertical arrows are surjective since all functors are exact.
	Hence, if $\theta^\pk_{a, X}$ is zero or surjective, so is $\theta^\pk_{a, Y}$.
	
	Consequently, $\theta^\pk_{a, Y}$ is zero if $a \not= \lambda_{R, \pk}+ \frac 1 2$, 
	and at least surjective when $a = \lambda_{R, \pk}+ \frac 12$.
	The injectivity in the latter case would follow from the fact that $\oplus_{a} \theta^\pk_{a, Y}$ is an injection given by compositions:
	\[\Gamma_{n, m-1} \circ T_{(\lambda_L, \lambda_R)}^{(\lambda_L, \lambda_R+ e_\pk)}(Y)
	\hookrightarrow
	\For \circ \Gamma_{n, m}(Y)
	\cong
	\bigoplus_{a\in \C} \Jac_a \circ \Gamma_{n, m}(Y).\]
	Now the conlusion follows.
\end{pf}

\subsection{On finite length objects}
	One can deduce Proposition \ref{nat-iso} for finite-length module from naturality of $\theta^\pk_a$ and five-lemma easily.
	
	In the meanwhile, any object in $\HC_\lambda$ is of finite length.
	Recall from Subsubsection \ref{subsubsec: Principle series} that irreducible objects in $\HC_\lambda$ are exhausted by finitely many $\bar X(w \lambda_L, \lambda_R)$, 
	where $w \in W(\g_0, \h_0) \cong S_n$ and $w\lambda_L - \lambda_R \in \C^n \cong \t^*$ is integral.
	If an object $X$ in $\HC_\lambda$ has infinitely many subquotients, then it must have one factor $\bar X(w \lambda_L, \lambda_R)$ repreating infinitely many times.
	Thus, $X$ has infinite dimensional $K$-isotypic subspaces, which contradicts to the admissibility in the definition of Harish-Chandra module.
	
	Hence, we have fully proved Proposition \ref{nat-iso} (assuming Lemma \ref{eigenvalue}).

\newcommand{\spann}{\operatorname{span}}
\section{Details in proof}\label{sec: eigenvalue}
	In this section we give a detailed proof of Lemma \ref{eigenvalue}.
	It is based mainly on the results in \cite{Chan23}.

\subsection{Filtration of $X(\lambda)\otimes V$}\label{Mackey}	
	Here we give the desired filtration in Lemma \ref{Std-under-Tran}.
	For this, we need the so-called Mackey isomorphism of parabolic induction,
	which leads to an isomorphism 
	$\Ind_B^G (\C_\lambda) \otimes V \cong \Ind_B^G( \C_\lambda \otimes V)$,
	$f \otimes z \to [g \mapsto f(g) \otimes z]$;
	see for example \cite[Theorem 2.103]{KV95}.
	
\begin{subequations}
	Consider the standard flag 
	\begin{equation}
	0 = V_0 \subset V_1 \subset \cdots \subset V_n = V = \C^n,
	\end{equation}
	where $V_k = \spann\{ e_1, \cdots, e_k\}$ consists of column vectors with first $k$ components non-zero only.
	This is also a filtration of $B$-representation, thus
	induces a filtration
	\begin{equation}
	 0 = X_0 \subset X_1 \subset \cdots \subset X_n = X,
	 \end{equation} 
	 of $X: = X(\lambda) \otimes V \cong \Ind_B^G (\C_\lambda \otimes V)$, where
	$X_k = \Ind_B^G( \C_\lambda \otimes V_k).$
\end{subequations}
	The parabolic induction is an exact functor, so
	\[X_k / X_{k-1} \cong \Ind_B^G \left( \C_\lambda \otimes (V_k / V_{k-1})\right).\]
	As a representation of $B$,
	$\C_\lambda \otimes (V_k / V_{k-1}) = \C_{(\lambda_L, \lambda_R + e_k )}$.
	Thus, as a representation of $G$, $X_k/ X_{k-1} \cong X(\lambda_L, \lambda_R+e_k)$.

\subsection{Generalized eigenvalues of $y_1$}\label{calculation-over-std}
	For the rest of Lemma \ref{eigenvalue}, we will verify $\forall F \in \Gamma_{n, m-1}(X_k)$,
\begin{equation}\label{action-of-y_1}
	y_1 \cdot F = (\lambda_{R, k} + \frac 1 2) F 
	\mod \Gamma_{n, m-1}(X_{k-1}).
\end{equation} 
	It follows that $y_1$ preserves each $\Gamma_{n, m-1}(X_k)$, and acts on $\Gamma_{n, m-1}( X_k/ X_{k-1})$ by the scalar $\lambda_{R, k} + \frac 1 2$.
	
	Since $G = BK$, $K \cap B = T$, there is an isomorphism of vector space
	\[
	\Phi: \Ind_B^G(\C_\lambda \otimes V) ^K \cong (\C_\mu \otimes V^{\otimes m})^T
	\]
	given by evaluation $F \mapsto F(1)$.
	As a $T$-representation, 
	\[V  = \bigoplus_{k=1}^n \C_{-e_k},\]
	where the weight space with respect to $-e_k \in \C^n \cong \t^*$ is spanned by $e_k \in \C^n = V$, thus
	$(\C_\mu \otimes V^{\otimes m})^T$ possesses a basis $\{1 \otimes e_\kappa\}$,
	where the index $\kappa$ runs over the maps $\{1, \cdots, m\} \to \{1, \cdots, n\}$ such that $\#\kappa^{-1}(k) = \mu_k$, $\forall k =1, \cdots, n$, and
	$e_\kappa =  e_{\kappa(1)} \otimes \cdots \otimes e_{\kappa(m)} \in V^{\otimes m}$.
	In particular, the space $\Gamma_{n, m}(X(\lambda))$ is non-zero only when each $\mu_k \geqslant 0$ and $\mu_1 + \cdots + \mu_n =m$;
	in this case it has dimension $\frac{m!}{\mu_1!\cdots \mu_n!}$.
	This is also calculated in \cite[Subsection 5.2]{Chan23}.
	In what follows \textbf{we will always assume} each $\mu_k = \lambda_{L, k} - \lambda_{R, k} \geqslant 0$.
	
	Let $F_\kappa \in \Gamma_{n, m}(X(\lambda))$ such that $\Phi(F_\kappa) = F_\kappa(1) = 1 \otimes e_\kappa$.
	It follows from definition that $\Gamma_{n, m-1}(X_k) =
	\Ind_B^G(\C_\lambda \otimes V_k \otimes V^{(m-1)})$ are sent to 
	$\left( \C_\mu \otimes V_k \otimes (V^{\otimes(m-1)}) \right)^T$ under $\Phi$, so
	this space is spanned by those $F_\kappa$ such that $\kappa(1) \leqslant k$.
	Thus, to verify \eqref{action-of-y_1}, it suffices to consider $F = F_\kappa$, and $k = \kappa(1)$;
	note that if $k > \kappa(1)$ then the equation holds trivially.
	We will divide these $F_\kappa$ into three parts:
\begin{subequations}
	\begin{align}
	F_\kappa	&:\kappa(1), \cdots \kappa(\mu_1) = 1,
	\label{eq: part-1}	\\
	F_\kappa	&:\kappa(1) = k>1, \kappa(2), \cdots, \kappa(\mu_1 + 1) = 1,
	\label{eq: part-2}	\\
	F_\kappa	&:\max\{\kappa(2), \cdots, \kappa(\mu_1)\} >1
	\label{eq: part-3}
	\end{align}
\end{subequations}
	Our proof will also be induction on $n$.
	Its starting point $n=1$ is easy. 
	In this case $X(\lambda) = \C_{(\lambda_L, \lambda_R)}$ is a one-dimensional representation of $\GL(1, \C) = \C^\times$,
	and $\Gamma_{n, m}(X(\lambda)) = \St(\Delta)$ is a one-dimensional module over $\H_m$, with $\Delta = [\lambda_R+ \frac 1 2, \lambda_L - \frac 1 2]$;
	this is \cite[Lemma 6.1]{Chan23}, and follows from direct computation by \eqref{action-of-H}.
	In particular, $y_1 \in \H_m$ acts by the scalar $\lambda_R + \frac 1 2$ on $\St(\Delta)$, as pointed out in Lemma \ref{Leb-of-Jac}.

\subsection{Part \eqref{eq: part-1}}
	This part follows directly from the results in \cite{Chan23}.
	Let 
	\begin{itemize}
	\item $m_1 = \mu_1$, $m_2 = \mu_2 + \cdots + \mu_n$,
	$V' = \spann_\C\{e_2, \cdots, e_n\}$, such that $V_1 \oplus V' = V$,
	\item $\mathcal X_1= \C e_1^{\otimes m_1} \otimes V'^{\otimes m_2} \subseteq V^{\otimes m}$, 
	$\mathcal W_1 = \Phi^{-1}(\C_\mu \otimes \mathcal X_1)^T$,
	\item $Y_1 = \C_{(\lambda_{L, 1}, \lambda_{R, 1})}$, $Y_2 = X(\lambda_L', \lambda_R') $, where
	$\lambda_L' = ( \lambda_{L, 2}, \cdots, \lambda_{L,n})$,
	$\lambda_R' = ( \lambda_{R, 2}, \cdots, \lambda_{R, n})$,
	\end{itemize}
	then $\mathcal W_1 \cong \Gamma_{1, m_1}(Y_1) \boxtimes \Gamma_{n-1, m_2}(Y_2)$ as $\H_{m_1} \otimes \H_{m_2}$-modules;
	this is Claim 2 in the proof of \cite[Theorem 5.7]{Chan23},
	since our definition of $\mathcal W_1$ is a special case of their $\mathcal W^{m_1, m_2}_1$.
	We already know $y_1 \in \H_{m_1}$ acts on
	$\Gamma_{1, m_1} (\C_{(\lambda_{L, 1} , \lambda_{R, 1})}) =
	\St\left( [\lambda_{R, 1} + \frac 1 2, \lambda_{L, 1} - \frac 1 2]\right)$
	by the scalar $\lambda_{R, 1} + \frac 1 2$.
	Then $y_1 \in \H_{m_1} \subseteq \H_m$ acts on 
	$\mathcal W_1 \cong \Gamma_{1, m_1}(Y_1) \boxtimes \Gamma_{n-1, m_2}(Y_2)$ 
	by the same scalar $\lambda_{R, 1} + \frac 1 2$.
	The subspace $\mathcal W_1 \subseteq \Gamma_{n, m}( X(\lambda))$ is spanned exactly by those $F_\kappa$ in Part \eqref{eq: part-1}.
	Thus, we obtain \eqref{action-of-y_1} for these $F_\kappa$.

\subsection{Part \eqref{eq: part-2}}
	Regard $S_m$ as the permutation group of the index set $\{1, \cdots, m\}$, then 
	according to \eqref{rmk: action-of-s}, $w\in S_m \subset \H_m$ acts on 
	$x \otimes v_1 \otimes \cdots \otimes v_m \in X(\lambda) \otimes V^{\otimes m}$ by 
	$(-1)^{l(w)} x \otimes v_{w^{-1}(1)} \otimes \cdots \otimes v_{w^{-1}(n)}$.
	For $F_\kappa$ in Part \eqref{eq: part-2},
	take $w = s_{m_1} \cdots s_1 \in S_m$ such that $w^{-1}(m_1+1)=1$, 
	$w^{-1}\{1, \cdots, m_1\} = \{2, \cdots, m_1+1\}$,
	then $w \cdot F_\kappa = (-1)^{m_1} F_{\kappa \circ w^{-1}}$ with $F_{\kappa\circ w^{-1}}$ lying in Part \eqref{eq: part-1}.
	
	It follows from induction hypothesis on $n$ that
	\begin{equation}\label{eq: induction}
	y_{m_1+1} \cdot F_{\kappa \circ w^{-1}}(1) = 
	\lambda_{R, k+1} F_{\kappa \circ w^{-1}}(1) \mod
	(\C_\mu \otimes \C e_1^{\otimes m_1} \otimes V_{k-1}' \otimes {V'}^{\otimes (m_2 -1)})^T,
	\end{equation}
	where $V_{k-1}' = V' \cap V_{k-1} = \spann_\C\{e_2, \cdots, e_{k-1}\}$, since
	\begin{itemize}
	\item $F_{\kappa \circ w^{-1}} \in \mathcal W_1$, and 
	$\mathcal W_1 \cong 
	\Gamma_{1, m_1}( \C_{(\lambda_{L, 1}, \lambda_{R, 1})}) \boxtimes
	\Gamma_{n-1, m_2}( X(\lambda_L', \lambda_R'))$
	as $\H_{m_1} \otimes \H_{m_2}$-modules,
	\item $X(\lambda_L', \lambda_R') = 
	\Ind_{B_2}^{G_2}(\C_{(\lambda_L', \lambda_R')})$,
	where $G_2 = \GL(n-1, \C)$, and $B_2 \subseteq G_2$ consists of upper triangluar matrices,
	\item $G_2$ acts on $V' \subset V$ as its conjugate standard representation,
	\item $\Gamma_{n-1, m_2}(X(\lambda_L', \lambda_R')) = 
	\Ind_{B_2}^{G_2}(
		\C_{(\lambda_L', \lambda_R')} \otimes {V'}^{\otimes m_2}
	)^{K_2}$ admits a filtration with terms
	\[\Ind_{B_2}^{G_2}(
		\C_{(\lambda_L', \lambda_R')} \otimes V'_k \otimes {V'}^{\otimes (m_2-1)}
	)^{K_2}, \quad
	k=1, \cdots, n,\]
	where $K_2 = \U(n-1) \subset \GL(n-1)$.
	\end{itemize}

	Following the defining relation of $\H_m$ given in Subsection \ref{subsec: H_m},
\[\begin{aligned}
	y_1 \cdot (-1)^{m_1}F_\kappa =&
	 y_1 \cdot s_1 \cdots s_{m_1} \cdot F_{\kappa\circ w^{-1}}
	=
	s_1 y_2 s_2 \cdots s_{m_1} \cdot F_{k} + s_2 \cdots s_{m_1} \cdot F_{\kappa\circ w^{-1}}\\
	=&
	s_1 \cdots s_{m_1} \cdot y_{m_1+1} \cdot F_{\kappa\circ w^{-1}} 
	+ \sum_{i=1}^{m_1} s_1 \cdots \hat s_i \cdots s_{m_1} \cdot F_{\kappa\circ w^{-1}}.
\end{aligned}\]
	It lies in $( \lambda_{R, k} + \frac 1 2) \cdot (-1)^{m_1} F_\kappa +
	\Ind_B^G( \C_\lambda \otimes V_{k-1} \otimes V^{\otimes (m-1)})^K$
	due to
	\begin{itemize}
	\item the equation \eqref{eq: induction},
	\item $s_1 \cdots s_{m_1} \cdot (\C_\mu \otimes \C e_1^{\otimes m_1} \otimes V_{k-1}' \otimes {V'}^{\otimes (m_2 -1)})^T =
	(\C_\mu \otimes V_{k-1}' \otimes  \C e_1^{\otimes m_1} \otimes {V'}^{\otimes (m_2 -1)})^T$, and
	\item $s_1 \cdots \hat s_i \cdots s_{m_1} \cdot F_{\kappa \circ w^{-1}}(1) \in
	s_1 \cdots \hat s_i \cdots s_{m_1} \cdot (\C_\mu \otimes \C e_1^{\otimes m_1} \otimes V_{k-1}' \otimes {V'}^{\otimes (m_2 -1)})^T \subseteq
	(\C_\mu \otimes \C e_1 \otimes {V}^{\otimes (m -1)})^T$.
	\end{itemize}
	Thus, we obtain \eqref{action-of-y_1} for these $F_\kappa$.

\subsection{Part \eqref{eq: part-3}}
	Let $S_{m-1} \to S_m$ be the subgroup fixing the first index $1$.
	If $F_\kappa$ belongs to Part \eqref{eq: part-3}, then $\kappa^{-1} (1)$ has at least $m_1-1$ common elements with $\{2, \cdots, m\}$.
	Take a $w \in S_{m-1}$ sending $\kappa^{-1}(1) \cap \{2, \cdots, m\}$ into $\{2, \cdots, m_1+1\}$.
	Thus, $w\cdot F_\kappa = (-1)^{l(w)} F_{\kappa \circ w^{-1}}$, where
	$F_{\kappa \circ w^{-1}}$ lies in Part \eqref{eq: part-1} or \eqref{eq: part-2}, and has been proved to satisfy \eqref{action-of-y_1}.
	Since any $w\in S_{m-1} \subset \H_m$
	\begin{itemize}
	\item commutes with $y_1\in \H_m$, and
	\item preserves the subspace
	$\Gamma_{n, m-1}(X_k) = (X_k \otimes V^{\otimes(m-1)})^K =
	(X(\lambda) \otimes V_k \otimes V^{\otimes(m-1)})^K$,
	\end{itemize}
	we deduce those $F_\kappa$ in Part \eqref{eq: part-3} also satisfies \eqref{action-of-y_1}.

\newcommand{\pr}{\operatorname{pr}}
\newcommand{\n}{{\mathfrak n}}
\section{A dual version}\label{dual}
	We now discuss the dual version of diagram \eqref{diagram: compability}:
	\begin{equation}\label{diagram: dual-compability}
	\begin{tikzcd}
	{\HC_{(\lambda_L, \lambda_R)}} \arrow[r, "{\Gamma_{n, m}}"] \arrow[d, "{T_{(\lambda_L, \lambda_R)}^{(\lambda_L - e_\pk, \lambda_R )}}"'] & 
	\H_m \Mod \arrow[d, "{\Jac^{\lambda_{L, \pk} - \frac 1 2}}"] \\
	{\HC_{(\lambda_L- e_\pk, \lambda_R)}} \arrow[r, "{\Gamma_{n, m-1}}"]                                                                    &
	\H_{m-1} \Mod                                             
	\end{tikzcd},
	\end{equation}
	It can be deduced from the following commutative diagrams:
	\[\begin{tikzcd}
	\HC_{(\lambda_L, \lambda_R)} \arrow[rd, "\sim" description] 
	\arrow[ddd] \arrow[rrr]	&                                                                 
	& 
	{}& 
	\H_m \Mod \arrow[ddd] \\
	& 
	\HC_{(-\overline{\lambda_L}, -\overline{\lambda_R})} 
	\arrow[r] \arrow[d] \arrow[ru, "c" description, Rightarrow]	& 
	\H_m \Mod \arrow[d] \arrow[ru, "\sim" description]	&                       
		\\
	& 
	\HC_{(-\overline{\lambda_R}, -\overline{\lambda_L} + e_\pk )} \arrow[r] \arrow[ru, "\theta" description, Rightarrow]      & 
	\H_{m-1} \Mod \arrow[rd, "\sim"] \arrow[ruu, "h" description, Rightarrow] &
	\\
	\HC_{(\lambda_L- e_\pk, \lambda_R)} \arrow[ru, "\sim" description] \arrow[rrr] \arrow[ruu, "r" description, Rightarrow] & 
	{} \arrow[ru, "c" description, Rightarrow]                      &                                                                          
	& 
	\H_{m-1} \Mod        
	\end{tikzcd}\]
	Here, each of the four $\sim$ is an equivalence of categories induced from Hermitian dual.
	The $\theta$ denotes the natural isomorphism for \eqref{diagram: compability},
	The two $c$ assert the compatibility between $\Gamma_{n, m}$ and Hermitian dual, as established in \cite[Theorem 7.5]{Chan23}.
	To complete the picture, it remains to explain the compatibility of Hermitian dual with translation and Jacquet functor,
	which are given by the natural isomorphism $r$ and $h$ respectively.

\subsection{On the real side}
	We follow the definition and notations in \cite[subsection 7.1]{Chan23}.
	Recall that the complex conjugate over $\g = \g_0 \otimes \C \cong \g_0 \times \g_0$ is $(E, E') \mapsto (\bar E', \bar E)$. 
	
	Define the \textbf{Hermitial dual} $X^h$ of a Harish-Chandra module $X$ of $\GL(n, \C)$ as follows.
	As a space, it consists of $K$-finite conjugate linear functionals over $X$, 
	with a non-degenerate sesquilinear pairing $\la-, -\ra: X^h \times X \to \C$;
	its structure of $(\g, K)$-module is determined by
	\[\la(E, E') \cdot f, x\ra = -\la f, (\bar E', \bar E) \cdot x\ra, \quad
	\la k \cdot f, x\ra = \la f, k^{-1} \cdot x\ra.\]
	An alternative definition of the Hermitian dual, which can be found in \cite[section VI.2]{KV95}, arises from the composition of contragredient functor and conjugation functor.
	It follows that the functor $X\mapsto X^h$ is exact, and $X^{hh} \cong X$ naturally.
	
	Inspired from \cite[(11.199)]{KV95} that 
	$X(\lambda_L, \lambda_R)^h \cong X(-\overline \lambda_R, - \overline \lambda_L),$
	one may infer the following effect of Hermitian dual on infinitesimal character.
\begin{lem}\label{Hermi-dual-and-inf-char}
	If the Harish-Chandra module $X$ has generalized infinitesimal character $(\lambda_L, \lambda_R)$, then its Hermitian dual $X^h$ has generalized infinitesimal character $(-\overline \lambda_R, -\overline \lambda_L)$.
\end{lem}
\begin{pf}
	Since $X \mapsto X^h$ is an exact fuctor, it sufficies to prove for the modules with an infinitesimal character.
	Denote the complex conjugate over $\g$ by $\sigma$.
	There is a conjugate-linear anti-isomorphism of the universal enveloping algebra $U(\g)$:
	$E_1\otimes \cdots \otimes E_r \mapsto (-\sigma)(E_r) \otimes \cdots \otimes (-\sigma)(E_1)$.
	Denote it by $-\sigma$, so $\forall u \in U(\g)$, and $f\in X^h$, $x\in X$, we have
	$\la u \cdot f, x\ra = \la f, (-\sigma)(u) \cdot x \ra$.
	
	Suppose the algebra $Z(\g)$ acts on $X$ via the multiplicative fuctional $\chi$, then it acts on $X^h$ via $\chi^h : = \overline{\chi \circ (-\sigma)}$.
	Here $-\sigma$ becomes a conjugate-linear isomorphism over the algebra $Z(\g)$.
	Alternatively speaking, we have commutativity of the first box in following diagrams:
	\[\begin{tikzcd}
	\C \arrow[d, "\textrm{conjugation}"'] & 
	Z(\g) \arrow[l, "\chi^h"'] \arrow[d, "-\sigma"] \arrow[r, "\pr_\n"] & 
	{\C[\h]} \arrow[d, "-\sigma"] \arrow[r, "T_\n"] & 
	{\C[\h]} \arrow[d, "-\sigma"] \\
	\C                                    & 
	Z(\g) \arrow[l, "\chi"'] \arrow[r, "\pr_{\n_-}"]                    & 
	{\C[\h]} \arrow[r, "T_{\n_-}"]                  & 
	{\C[\h]}                     
	\end{tikzcd}.\]
	Let's explain the rest of this diagram.
	The algebra isomorphism $-\sigma$ over $\C[\h]$ is also induced from  $-\sigma$ over $U(\g)$.
	The algebra homomorphisms $\pr_{\n}$ and $T_{\n}$ compose into the Harish-Chandra isomorphism $Z(\g) \cong \C[\h]^W$.
	More precisely, choose ($\sigma$-invariant) nilpotent subalgebras $\n, \n_-$ such that $\h \oplus \n$, $\h \oplus \n_-$ are opposite Borel subalgebras of $\g$;
	they exist since $\g = \g_0 \otimes_\R \C$ and $\g_0$ has complex structure.
	Define 
\begin{itemize}
	\item 
	$\pr_\n(z)$ for $z\in Z(\g)$ to be the unique element in $\C[\h]$ such that $z \equiv \pr_\n(z) \mod U(\g) \n$, or equivalently, 
	 the unique element in $\C[\h]$ such that $z \equiv \pr_\n(z) \mod \n_-U(\g)$, see \cite[Proposition 5.34 and (5.37)]{Knapp96}, and	 
	\item 
	$T_\n(H)$ for $H \in \h$ to be $H - \rho_\n(H) \in \C[\h]$,
	where $\rho_\n$ is half the sum of positive roots with respect to the Borel subalgebra $\h \oplus \n$,
	which natrually induces an algebra isomorphism over $\C[\h]$.
\end{itemize}
	There are similar definitions for $\pr_{\n_-}$ and $T_{\n_-}$.
	\cite[Theorem 5.44]{Knapp96} shows that both $T_\n \circ \pr_n$ and $T_{\n_-} \circ \pr_{\n_-}$ establish the same isomorphism between $Z(\g)$ and $\C[\h]^W$.
	It is exactly the Harish-Chandra isomorphism.
	The commutativity of the above diagram follows from easy calculation.
	
	Consequently, the infinitesimal character $\lambda^h$ of $X^h$, as a (Weyl group orbit of) wieght in $\h^*$, would be related to $\lambda$ via the equation:
	$\lambda^h = \overline{\lambda \circ (-\sigma)}.$
	More explicitly, let $(H, H') \in \h$, and $\lambda = (\lambda_L, \lambda_R)$, then
	\[\lambda^h(H, H') = \overline{\lambda (-\bar H', -\bar H)} 
	= \overline{-\lambda_L(\bar H') - \lambda_R(\bar H)}
	= -\overline \lambda_R(H) - \overline \lambda_L( H').\]
	Hence, we obtain $\lambda^h = ( -\overline \lambda_R, - \overline \lambda_L)$.	
\end{pf}
	
	As a consequence, the Hermitian dual sends $\HC_{(\lambda_L, \lambda_R)}$ to $\HC_{(-\overline \lambda_R, -\overline\lambda_L )}.$
	The natural isomorphism $r$ can be induced by the following diagram:
\[\begin{tikzcd}
	{\HC_{(\lambda_L, \lambda_R)}} && \HC_n && {\HC_{(\lambda_L, \lambda_R+e_\pk)}} \\
	\\
	{\HC_{(-\overline\lambda_R, -\overline\lambda_L)}} && \HC_n && {\HC_{(-\overline\lambda_R-e_\pk, -\overline\lambda_L)}}
	\arrow["{-\otimes V}", from=1-1, to=1-3]
	\arrow["\sim"', from=1-1, to=3-1]
	\arrow["{P_{(\lambda_L, \lambda_R+e_\pk)}}", from=1-3, to=1-5]
	\arrow["\sim", from=1-3, to=3-3]
	\arrow["\sim", from=1-5, to=3-5]
	\arrow[Rightarrow, from=3-1, to=1-3]
	\arrow["{-\otimes V^h}"', from=3-1, to=3-3]
	\arrow[Rightarrow, from=3-3, to=1-5]
	\arrow["{P_{(-\overline\lambda_R-e_\pk, -\overline\lambda_L)}}"', from=3-3, to=3-5]
\end{tikzcd}.\]
	Note that the Hermitian dual of conjugate standard representation $V$ is given by $V^h = \C^n$, with $g\in G$ acting by left multiplication of $(g^t)^{-1}$, and
	has weights $(-e_i, 0) \in \C^n \oplus \C^n \cong \h_0^* \oplus \h_0^*$.

\subsection{On the $p$-adic side}
	We follow the definition and notations in \cite[subsection 7.2]{Chan23}.
	There is a conjugate linear anti-involution $^*$ over $\H_m$, determined by
	$s_i^*=s_i$, $y_i^*=-w_0(y_{m+1-i})w_0^{-1},$
	where $w_0 \in S_m$ is the longest element in it.

	Define the \textbf{Hermitian dual} $M^*$ of  an $\H_m$-module $M$ as follows.
	As a space, it consists of conjugate linear functionals from $M$ to $\C$; 
	its structure of $\H_m$-module is determined by
	\[(h\cdot f)(x)= f(h^*\cdot x).\]
	In particular, there is a non-degenerate sesquilinear pairing $\la -, - \ra: M^*\times M \to \C$ such that
	$\la h\cdot f, x \ra=\la f, h^*\cdot x \ra$.
	

	Furthermore, the following lemma is well-known, and can be verified by direct computation:
\begin{lem}\label{Hermi-dual-and-Jac}
	There is a natrual isomorphism fitting into the following diagram of functors:
\[\begin{tikzcd}
	{\H_m\Mod} && {\H_m\Mod} \\
	\\
	{\H_{m-1}\Mod} && {\H_{m-1}\Mod}
	\arrow["\sim", from=1-1, to=1-3]
	\arrow["{\Jac_a}"', from=1-1, to=3-1]
	\arrow["{\Jac^{-\bar a}}", from=1-3, to=3-3]
	\arrow[Rightarrow, from=3-1, to=1-3]
	\arrow["\sim"', from=3-1, to=3-3]
\end{tikzcd}.\]
	Informally speaking, the Hermitian dual of generalized $a$-eigenspace with respecct to $y_1$ is natrually isomorphic to the generalized $(-\bar a)$-eigenspace of Hermitian dual with respect to $y_m$.
\end{lem}	
\begin{pf}
	Two ways of embedding $\H_1 \otimes \H_{m-1}$ into $\H_m$ are involved in this lemma.
	To be precise, denote the polynomial generators and group generators of $\H_{m-1}$ by $z_i, t_i$ (instead of $y_i, s_i$) respectively,
	and let $\H_1 = \C[z]$.
	The two embeddings of $\H_1 \otimes \H_{m-1}$ into $\H_m$ are given by
	\[
	\iota_1: z\mapsto y_1, z_i \mapsto y_{i+1}, t_i \mapsto s_{i+1};\quad
	\iota_m: z\mapsto y_m, z_i \mapsto y_{i}, t_i \mapsto s_{i}.	
	\]
	Then one can check the following by definition:
\begin{itemize}
	\item $\Jac_a$ is given by restricting $\H_m$-module to $\H_1\otimes \H_{m-1}$-module via $\iota_1$, 
	and then taking generalized $a$-eigenspace with respecct to $z \in \H_1$;
	\item $\Jac^{-\bar a}$ is given by restricting $\H_m$-module to $\H_1\otimes \H_{m-1}$-module via $\iota_m$, 
	and then taking generalized $(- \bar a)$-eigenspace with respecct to $z \in \H_1$.
\end{itemize}
	
	Let $M$ be an $\H_m$-module, and $M^*$ be its Hermitial dual. 
	By definition, there is a non-degenerate sesquilinear pairing $\la -, - \ra_m: M^*\times M \to \C$ such that $\forall h \in \H_m$, 
	$\la h\cdot f, x \ra_m=\la f, h^*\cdot x \ra_m$.
	Restricting $M$ and $M^*$ to $\H_1 \otimes \H_{m-1}$-modules via $\iota_1$ and $\iota_m$ respectively.
	We would like to define a non-degenerate sesquilinear pairing $\la -, - \ra_{m-1}: M^*\times M \to \C$ such that $\forall h \in \H_1\otimes \H_{m-1}$, 
	$\la \iota_m(h)\cdot f, x \ra_{m-1}=\la f, \iota_1(h^*)\cdot x \ra_{m-1}$,
	where the conjugate linear anti-involution $^*$ over $\H_1\otimes \H_{m-1}$ is induced from those $^*$  over $\H_1$ and $\H_{m-1}$ natrually.
	
	Once such a $\la-, - \ra_{m-1}$ is defined, then by
	$\la \iota_m(z) \cdot f, x \ra_{m-1} = \la f, \iota_1(-z) \cdot x \ra_{m-1}$
	we obtain a non-degenerate sesquilinear pairing between generalized $a$-eigenspace of $M$ and  generalized $(- \bar a)$-eigenspace of $M^*$ with respecct to $z$.
	Moreover, by
	$\la \iota_m(h) \cdot f, x \ra_{m-1} = \la f, \iota_1(h^*) \cdot x \ra_{m-1}$
	for $h\in \H_{m-1}$ we know such two subspaces becomes Hermitian dual of each other as $\H_{m-1}$-modules.
	This is exactly what we want.
	
	Denote the longest element in $S_{m-1}$ by $u_0$, and define
	\[\la f, x \ra_{m-1} := \la f, w_0 \iota_1(u_0^{-1}) \cdot x \ra_m.\]
	The verification of $\la \iota_m(h) \cdot f, x \ra_{m-1}=\la f, \iota_1(h^*) \cdot x \ra_{m-1}$ are left to reader.	
\end{pf}
\begin{rmk}
	Here we scketch another, perhaps easier proof.
	Let $(\pi, M)$ be an $\H_m$-module.
	Its inner twist $\pi^w$ by $w \in S_m$ belongs to the same isomorphic class of $\pi$;
	here $\pi^w(h) := \pi( wh w^{-1})$, and
	the intertwining operator $\pi \stackrel \sim \to \pi^w$ can be given by $x \mapsto \pi(w) x$.
	As a consequence, we can use another conjugate linear anti-involution $^\star$ over $\H_m$ for an equivalent definition of Hermitian dual,
	which would simplify the verification of Lemma \ref{Hermi-dual-and-Jac}.
	
	Take $s_i^\star= s_{m-i},$ $y_i^\star=- y_{m+1-i}$, and one calculates that $\forall h \in \H_1 \otimes \H_{m-1}$, $\iota_1(h)^\star = \iota_m(h^\star)$.
	Then the Hermitial dual $M^\star$ of $M$ as $\H_m$-module is also its Hermitial dual as $\H_1 \otimes \H_{m-1}$-module,
	so are their generalized eigenspaces with respect to $z\in \H_1$ respectively.
\end{rmk}

	This lemma provides the desired natrual isomorphism $h$.

\newcommand{\Perv}{{\mathrm{Perv}}}
\newcommand{\IC}{{\mathrm{IC}}}
\newcommand{\Rep}{{\mathrm{Rep}}}

\newcommand{\FDR}{{\mathcal S}}
\newcommand{\Speh}{{\mathrm{Speh}}}
\newcommand{\Tran}{{\mathrm T}}
\newcommand{\Herm}{{\mathrm H}}
\newcommand{\ES}{{\mathcal {E}}}
\newcommand{\MRV}{{\mathcal D}}
\newcommand{\AV}{{\mathrm{AV}}}
\newcommand{\bp}{\mathrm{bp}}
\newcommand{\Std}{{\mathrm {Std}}}
\newcommand{\q}{{\mathfrak q}}
\newcommand{\p}{{\mathfrak p}}
\renewcommand{\l}{{\mathfrak l}}
\renewcommand{\k}{{\mathfrak k}}
\renewcommand{\u}{{\mathfrak u}}
\renewcommand{\b}{{\mathfrak b}}
\newcommand{\z}{{\mathfrak z}}
\renewcommand{\r}{{\mathfrak r}}
\renewcommand{\v}{{\mathfrak v}}

\bibliographystyle{alpha}
\bibliography{D:/LaTeX/Research/Reference/LLC_p_adic.bib, 
	D:/LaTeX/Research/Reference/A_packet_construction_p_adic.bib, 
	D:/LaTeX/Research/Reference/A_packet_construction_archimedean.bib, 
	D:/LaTeX/Research/Reference/Real_Lie_group.bib, 
	D:/LaTeX/Research/Reference/LLC_archimedean.bib,
	D:/LaTeX/Research/Reference/Personal,
	Reference.bib}

\vspace{1em}
\begin{flushleft} \small
	Chang Huang: Yau Mathematical Sciences Center, Tsinghua University, Haidian District, Beijing 100084, China. \\
	E-mail address: \href{mailto:hc21@mails.tsinghua.edu.cn}{\texttt{hc21@mails.tsinghua.edu.cn}}
\end{flushleft}
\end{document}